# STUDENT EXPLANATION IN MIDDLE AND SECONDARY MATHEMATICS AND STATISTICS: A SCOPING LITERATURE REVIEW

Huixin Gao, Tanya Evans, and Anna Fergusson

University of Auckland

*This scoping review examines the literature on student explanation strategies in middle and secondary mathematics and statistics education from 2014 to 2024. Following the PRISMA protocol, we analyzed 41 studies that met the inclusion criteria. The findings classify student explanations into two types: self-explanation (SE) and peer explanation (PE). Both approaches enhance conceptual understanding and procedural knowledge, though each offers distinct benefits and challenges. SE is particularly effective when combined with worked examples and varies with prompting strategies, while PE significantly impacts students' affective development and social learning. The review identifies a significant gap in studies comparing the effectiveness of SE and PE, alongside an almost complete absence of research within statistics education.*

## 1. INTRODUCTION

Over the past few decades, researchers have emphasized the importance of student-constructed explanations in mathematics and statistics education (Rittle-Johnson, 2024). Student explanation requires learners to articulate their understanding of mathematical concepts and reasoning processes (Rittle-Johnson & Loehr, 2017). Research shows this approach enhances both conceptual understanding and core competencies like logical reasoning and abstract thinking (Rittle-Johnson et al., 2017; Fiorella & Mayer, 2016a), reflecting modern mathematics education's emphasis on active learning practices (Fiorella & Mayer, 2016b).

Despite the recognized value of student explanations in mathematics learning, research lacks systematic organization regarding its effective implementation in educational settings. This study aims to systematically analyze existing literature to clarify the current state of research on this strategy.

## 2. THEORETICAL GROUNDING

Student explanations can be understood through multiple theoretical frameworks. Human cognitive architecture theory (Atkinson & Shiffrin, 1968) describes three key stages: sensory memory for receiving inputs, working memory for active processing, and long-term memory for storing information in schemas. Related theories, such as cognitive load theory (Sweller, 1988), further explain that working memory has limited capacity in processing new instructional content. Explanation generation supports schema construction by directing cognitive resources to integrate new knowledge with existing understanding (Mayer & Moreno, 2003). In mathematics and statistics learning, this integration occurs through connecting concepts and applying principles to new contexts, fostering meaningful knowledge construction through cognitive processing and memory integration (Fiorella & Mayer, 2014; Lachner et al., 2021).

Social Presence Theory (Short et al., 1976) adds a social interaction perspective, emphasizing how learners adjust their explanations based on their audience's knowledge level—a process that deepens learning (Monrose Mills et al., 2020). When implemented through peer teaching and group problem-solving, these activities enhance social presence, motivation, and cognitive engagement (Nasir et al., 2023; Morris et al., 2023; Alegre et al., 2019), leading to more effective memory formation.

These frameworks complement each other in explaining student explanations' effectiveness: cognitive architecture theory addresses information processing and knowledge construction, while social presence theory focuses on social interaction benefits. Together, they reveal how student explanations enhance learning through the synergy of cognitive processing and social interaction.

## 3. METHODS

This review follows Arksey & O'Malley's (2005) five-stage framework: (1) identifying the research question, (2) identifying relevant studies, (3) study selection, (4) charting the data, and (5) collating and reporting results. This framework ensures a rigorous and replicable review process.

### 3.1 Identifying the initial research question

Our review focuses on examining the specific use of student explanations in middle and secondary mathematics and statistics education. Therefore, we propose the following research question: How can the literature on student explanation strategies in middle and secondary mathematics and statistics education be categorized to reveal current research patterns and trends?

### 3.2 Identifying relevant studies

All publications included in this review were published in English and indexed in three major academic databases: ERIC (Education Resources Information Center), Scopus, and Web of Science. To ensure comprehensive coverage, we searched major databases, supplemented with Google Scholar and Research Rabbit searches, and conducted forward and backward citation tracking. The search used systematic keywords from preliminary review and expert consultation, as shown in Table 1.

We reviewed literature from 2014-2024, focusing on "explanation" and related terms in middle and secondary education (Grades 6-13). While our search might have missed studies using terms like "articulation" or "justification," we mitigated this through broad initial screening. We included only English-language, peer-reviewed articles.

| Search string |
|---|
| ("Secondary Education" OR "Middle School" OR "Intermediate School" OR "High School" OR "Junior High" OR "Senior High" OR "Grade 7" OR "Grade 8" OR "Grade 9" OR "Grade 10" OR "Grade 11" OR "Grade 12" OR "Adolescents" OR "Teenagers") AND |
| ("Mathematics" OR "Statistics" OR "STEM") AND |
| ("Learning by explaining" OR "Explaining-based learning" OR "Explanation-driven learning" OR "Self-explanation" OR "Self-explaining" OR "Learning by teaching" OR "Learning-by-teaching" OR "Team-based-explanation" OR "Peer tutoring" OR "Peer instruction" OR "Peer assist" OR "Peer-led") |

Table 1: Search string

### 3.3 Study selection

Following PRISMA protocol (Moher et al., 2009), our initial search retrieved 231 publications. After removing duplicates and screening titles and abstracts, 57 publications proceeded to full-text assessment. Further evaluation excluded 16 publications, leaving 41 studies for qualitative analysis (complete literature list is available on figshare: https://doi.org/10.17608/k6.auckland.26933131).

## 4. FINDINGS AND DISCUSSION

Based on the analysis and synthesis of the results of selected literature, the student explanation strategies can be categorized into two broad types: Self-explanation (SE) and Peer explanation (PE). Of the 41 studies reviewed, 16 focus on SE (16 on mathematics, 0 on statistics), while 25 discuss PE (24 on mathematics, and 1 on statistics).

### 4.1 Self-explanation

Self-explanation involves learners generating explanations to themselves (Hodds et al., 2014), actively connecting new information with existing knowledge (Rittle-Johnson et al., 2017). Students articulate their understanding of concepts, procedures, and problem-solving steps, helping identify knowledge gaps and build coherent mental models (Lachner et al., 2021; Rittle-Johnson & Loehr, 2017). In middle and secondary education, SE is implemented through written explanations in textbooks or student sheets, or through oral explanations via audio recordings or self-narration.

Six of the 16 studies examined combining SE with worked examples (e.g., Barbieri et al., 2021). This combination reduces errors, speeds up problem-solving, and decreases the need for teacher assistance (e.g., McGinn et al., 2015), while helping students overcome insufficient prior knowledge in algebra (e.g., Barbieri et al., 2023). However, students may become overly dependent on worked examples, potentially limiting knowledge transfer (Schalk et al., 2018).

Student achievement through SE manifests primarily in conceptual and procedural knowledge development (e.g., Özcan, 2024). For example, in mathematics learning, when students used diagrammatic SE as an anticipatory approach, they improved their procedural knowledge by better applying formal algebraic problem-solving strategies to transfer problems with negative numbers (Nagashima et al., 2021). On the conceptual side, SE through GeoGebra enhanced students' ability to visualize and connect different mathematical representations, demonstrating deeper conceptual understanding (Nordlander, 2022). Notably, females seem to benefit more from SE than males (Nguyen et al., 2022). This gender difference appears to stem from different engagement patterns: females tend to produce more thorough and thoughtful explanations, while males often provide quick, superficial responses, potentially missing the learning benefits that come from struggling with and articulating explanations (McLaren et al., 2022)

The self-explanation principle suggests that SE falls along a continuum: open-ended (students are simply prompted to self-explain, e.g., "Explain how you solved this problem") , focused (guiding students toward specific aspects, e.g., "Explain why you chose this mathematical operation"), scaffolded (providing support for explanation, e.g., "Using the given formula, explain each step of your solution"), resource-based (referring to specific materials, e.g., "Using the diagram provided, explain how you arrived at your answer"), and menu-based (selecting from predetermined explanations, e.g., "Choose from multiple explanation options") (McLaren et al., 2022). Menu-based SE, which guides thinking through preset options, is more helpful for students with weak foundations (Wong et al., 2019). However, in terms of long-term effects, open-ended SE, which allows students to freely establish knowledge connections, shows better learning outcomes (McLaren et al., 2022).

However, SE has its limitations. In complex tasks, excessive SE may increase cognitive load (Hänze & Leiss, 2022) and may even lead to cognitive confusion and pseudo-confidence in students, affecting effective knowledge construction (Andres et al., 2023).

**4.2 Peer explanation**

Peer explanation (PE) has emerged as a key strategy in mathematics education's shift toward active learning (Alegre et al., 2019b). In PE, students work in pairs to explain mathematical concepts and problem-solving strategies (Moliner & Alegre, 2020b, 2022). Unlike simple peer tutoring, PE involves reciprocal knowledge construction where students alternate between explainer and listener roles, collaboratively building mathematical understanding (Moliner & Alegre, 2020a). This benefits both parties: explainers deepen understanding through articulation, while listeners gain new perspectives through engagement (Martí Arnándiz et al., 2022).

The 25 PE studies reveal two main positive impacts. First, PE improves students' affective domain, enhancing self-concept and reducing mathematics anxiety (e.g., Alegre & Moliner, 2017) by creating a safe, inclusive environment that encourages participation (Arthur et al., 2022). Second, PE develops conceptual and procedural knowledge (Fukuda & Manalo, 2024) as students explain both how and why solutions work, motivating them to become role models (Roberts & Spangenberg, 2020). Same-age PE proves more effective than cross-age PE (Moliner & Alegre, 2020b, 2022), with female and younger students showing higher gains in self-efficacy (Arnándiz et al., 2022).

However, PE has several limitations: First, learning through teaching may lead students to spend considerable time on irrelevant activities, creating extraneous cognitive load and thus affecting learning outcomes (Fiorella et al., 2019). Second, PE may negatively impact explainers, leading them to develop a fixed mindset about intelligence and affecting their self-efficacy (Gandolfi et al., 2024).

Furthermore, PE's effectiveness is constrained by students' knowledge base and explanatory abilities, where weaker students' content knowledge may lead to misconceptions during transmission (Thomas et al., 2015). Additionally, PE's effectiveness largely depends on the paired students' knowledge levels. When knowledge differences between pairs are either too large or too small, PE's effectiveness may be compromised, with weaker students showing lower self-concept and participation willingness, while students with similar knowledge levels may benefit less from PE (Moliner & Alegre, 2020b).

## 5. CONCLUSION

This scoping review underscores the well-documented positive effects of self-explanation and peer explanation in mathematics education, where both strategies have been shown to enhance conceptual understanding and problem-solving skills. However, research on student explanations in statistics education is notably limited, with only one of 41 studies addressing PE in statistics learning, highlighting a significant research gap at the middle and secondary levels.

While self-explanation has been extensively studied for its cognitive benefits and peer explanation for its social and collaborative advantages, few studies have directly compared their relative impacts on learning outcomes. Future research should focus on comparative studies to identify the conditions and learner characteristics that influence the effectiveness of self-explanation and peer explanation. Such efforts are crucial to inform evidence-based teaching practices that optimize the use of these strategies to support student learning in both mathematics and statistics education.

## References

Alegre, F., & Moliner, L. (2017). Emotional and cognitive effects of peer tutoring among secondary school mathematics students. *International Journal of Mathematical Education in Science and Technology*, *48*(8), 1185–1205.

Alegre, F., Moliner, L., Maroto, A., & Lorenzo-Valentin, G. (2019b). Peer tutoring in algebra: A study in Middle school. *Journal of Educational Research*, *112*(6), 693–699.

Andres, A., Cloude, E. B., Baker, R., Alexandra Andres, J. M., Baker, R. S., & Lee, S. (2023). Investigating Cognitive Biases in Self-Explanation Behaviors during Game-based Learning about Mathematics. *Proceedings of the 31st International Conference on Computers in Education. Asia-Pacific Society for Computers in Education*.

Arksey, H., & O'Malley, L. (2005). Scoping studies: towards a methodological framework. *International Journal of Social Research Methodology*, *8*(1), 19–32.

Arthur, Y. D., Boadu, S. K., & Asare, B. (2022). Effects of Peer Tutoring, Teaching Quality and Motivation on Mathematics Achievement in Senior High Schools. *International Journal of Educational Sciences*, *37*(1–3), 35–43.

Atkinson, R. C., & Shiffrin, R. M. (1968). Human memory: A proposed system and its control processes. In K. W. Spence & J. T. Spence (Eds.), *The psychology of learning and motivation* (Vol. 2, pp. 89-195). Academic Press.

Barbieri, C. A., Booth, J. L., Begolli, K. N., & McCann, N. (2021). The effect of worked examples on student learning and error anticipation in algebra. *Instructional Science*, *49*(4), 419–439.

Barbieri, C. A., Booth, J. L., & Chawla, K. (2023). Let's be rational: worked examples supplemented textbooks improve conceptual and fraction knowledge. *Educational Psychology*, *43*(1), 1–21.

Fiorella, L., Kuhlmann, S., & Vogel-Walcutt, J. J. (2019). Effects of Playing an Educational Math Game That Incorporates Learning by Teaching. *Journal of Educational Computing Research*, *57*(6), 1495–1512.

Fiorella, L., & Mayer, R. E. (2014). Role of expectations and explanations in learning by teaching. *Contemporary Educational Psychology*, *39*(2), 75–85.

Fiorella, L., & Mayer, R. E. (2016a). Effects of Observing the Instructor Draw Diagrams on Learning From Multimedia Messages. *Journal of Educational Psychology*, *108*(4), 528–546.

Fiorella, L., & Mayer, R. E. (2016). Eight ways to promote generative learning. *Educational Psychology Review*, *28*(4), 717-741.

Fukuda, M., & Manalo, E. (2024). Promoting learners' self-regulated textbook use for overcoming impasses in solving mathematics exercises. *Research in Mathematics Education*, *26*(1), 133–155.

Gandolfi, E., Limata, T., Favatà, R., & Ianì, F. (2024). Does peer tutoring have negative effects? An investigation and intervention on tutors' implicit theories and beliefs of intelligence. *Educational Psychology*, *44*(5), 632–648.

Hänze, M., & Leiss, D. (2022). Using heuristic worked examples to promote solving of reality-based tasks in mathematics in lower secondary school. *Instructional Science*, *50*(4), 529–549.

Lachner, A., Jacob, L., & Hoogerheide, V. (2021). Learning by writing explanations: Is explaining to a fictitious student more effective than self-explaining? *Learning and Instruction*, *74*.

Arnándiz, O. M., Moliner, L., & Alegre, F. (2022). When CLIL is for all: Improving learner motivation through peer-tutoring in Mathematics. *System*, *106*, 102773.

Mayer, R. E., & Moreno, R. (2003). Nine ways to reduce cognitive load in multimedia learning. *Educational Psychologist*, *38*(1), 43–52.

McGinn, K. M., Lange, K. E., & Booth, J. L. (2015). A Worked Example for Creating Worked Examples. *Mathematics Teaching in the Middle School*, *21*(1), 26–33.

McLaren, B. M., Richey, J. E., Nguyen, H. A., & Mogessie, M. (2022). A Digital Learning Game for Mathematics that Leads to Better Learning Outcomes for Female Students: Further Evidence. *Proceedings of the 16th European Conference on Games Based Learning*, *16*(1).

Moher, D., Liberati, A., Tetzlaff, J., Altman, D. G., & Group, and the P. (2009). Reprint—Preferred Reporting Items for Systematic Reviews and Meta-Analyses: The PRISMA Statement. *Physical Therapy*, *89*(9), 873–880.

Moliner, L., & Alegre, F. (2020a). Effects of peer tutoring on middle school students' mathematics self-concepts. *PLoS ONE*, *15*(4).

Moliner, L., & Alegre, F. (2020b). Peer Tutoring Effects on Students' Mathematics Anxiety: A Middle School Experience. *Frontiers in Psychology*, *11*.

Moliner, L., & Alegre, F. (2022). Peer tutoring in middle school mathematics: academic and psychological effects and moderators. *Educational Psychology*, *42*(8), 1027–1044.

Moliner, L., Alegre, F., & Lorenzo-Valentín, G. (2022). Peer Tutoring and Math Digital Tools: A Promising Combination in Middle School. *Mathematics*, *10*(13).

Monrose Mills, N., Blackmon, A., McKayle, C., Stolz, R., & Romano, S. (2020). *Peer-Led Team Learning and its effect on mathematics self-efficacy and anxiety in a developmental mathematics course* (pp. 93–102).

Morris, P., Agbonlahor, O., Winters, R., & Donelson, B. (2023). Self-efficacy curriculum and peer leader support in gateway college mathematics. *Learning Environments Research*, *26*(1), 219–240.

Nagashima, T., Bartel, A. N., Yadav, G., & Tseng, S. (2021, June). Using Anticipatory Diagrammatic Self-Explanation to Support Learning and Performance in Early Algebra. *International Society of the Learning Sciences*.

Nasir, R., Fadzli, N., Rusli, S. A. M., Kamaruzzaman, N. S., Sheng, N. N., Mohammad, V. Y. Z., & Shukeri, N. H. H. (2023). Peer tutoring learning strategies in mathematics subjects: Systematic literature review. *European Journal of Educational Research*, *12*(3), 1407–1423.

Nguyen, H. A., Hou, X., Richey, J. E., & McLaren, B. M. (2022). The Impact of Gender in Learning With Games. *International Journal of Game-Based Learning*, *12*(1), 1–29.

Nordlander, M. C. (2022). Lifting the understanding of trigonometric limits from procedural towards conceptual. *International Journal of Mathematical Education in Science and Technology*, *53*(11), 2973–2986.

Özcan, Z. Ç. (2024). Evaluating the Impact of the WEI4S Instructional Approach on Middle School Students' Algebraic Problem-Solving Skills. *Education Sciences*, *14*(1).

Rittle-Johnson, B., & Loehr, A. M. (2017). Eliciting explanations: Constraints on when self-explanation aids learning. *Psychonomic Bulletin and Review*, *24*(5), 1501–1510.

Rittle-Johnson, B., Loehr, A. M., & Durkin, K. (2017). Promoting self-explanation to improve mathematics learning: A meta-analysis and instructional design principles. *ZDM - Mathematics Education*, *49*(4), 599–611.

Rittle-Johnson, B. (2024). Encouraging students to explain their ideas when learning mathematics: A psychological perspective. *The Journal of Mathematical Behavior*, *76*, 101192.

Roberts, A. K., & Spangenberg, E. D. (2020). Peer tutors' views on their role in motivating learners to learn mathematics. *Pythagoras*, *41*(1), 1–13.

Schalk, L., Schumacher, R., Barth, A., & Stern, E. (2018). When problem-solving followed by instruction is superior to the traditional tell-and-practice sequence. *Journal of Educational Psychology*, *110*(4), 596–610.

Short, J., Williams, E., & Christie, B. (1976). *The social psychology of telecommunications*. John Wiley & Sons.

Sweller, J. (1988). Cognitive Load During Problem Solving: Effects on Learning. *Cognitive Science*, *12*(2), 257–285.

Thomas, A. S., Bonner, S. M., Everson, H. T., & Somers, J. A. (2015). Leveraging the power of peer-led learning: investigating effects on STEM performance in urban high schools. *Educational Research and Evaluation*, *21*(7–8), 537–557.

Wong, R. M., Adesope, O. O., & Carbonneau, K. J. (2019). Process-and Product-Oriented Worked Examples and Self-Explanations to Improve Learning Performance. *Journal of STEM Education: Innovations and Research*, *20*(2).